\documentclass[12pt]{amsart}
\usepackage{amssymb,mathrsfs,color,bbold}

\makeatletter

\newcommand{\Rmnum}[1]{\expandafter\@slowromancap\romannumeral #1@}
\makeatother
\usepackage{mathtools}

\theoremstyle{plain}

\theoremstyle{definition}

\usepackage[dvips]{graphics}
\usepackage{mathrsfs}
\usepackage{amssymb}
\usepackage{stmaryrd}
\usepackage{pifont}
\usepackage{amsmath}
\usepackage{hyperref}
\usepackage{amsthm}
\usepackage{cases}
\usepackage{mathtools}
\allowdisplaybreaks 
\def\beqn{\begin{displaymath}}
\def\eeqn{\end{displaymath}}

\def\bea{\begin{eqnarray}}
\def\eea{\end{eqnarray}}
\def\bean{\begin{eqnarray*}}
\def\eean{\end{eqnarray*}}

\def\beq{\begin{equation}}
\def\eeq{\end{equation}}

\begin{document}

\begin{center}
{\Large   $\sigma-$Unbounded Dunford-Pettis operators on Banach lattices}
\bigskip

\ Miao Ouyang, Zili Chen,Jingxi Chen, Zhangjun Wang*

School of Mathematics, Southwest Jiaotong University,
Chengdu, Sichuan,
China, 610000 \ \
\end{center}

\noindent\textbf{Abstract} \ An operator $T$ from a Banach lattice $E$ into a Banach lattice $F$ is called unbounded Dunford-Pettis operator whenever $x_{n}\xrightarrow{uaw}0$ implies
$Tx_n\xrightarrow{un}0$.
In this paper, we study some properties of the operators and relationships with the other classes of operators.

\textbf{2010 MR Subject Classification}
46A40, 46B42\\

\bigskip

\textbf{1. Introduction and preliminaries}
\bigskip

Several recent papers investigated unbounded convergence. A net $(x_\alpha)$ in Banach lattices $E$ is \emph{unbounded order (resp. norm, absolute weak) convergent} to $x\in E$ if $(|x_\alpha-x|\wedge u)$ is order (resp. norm, absolute weak) convergent to $x$ for any $u\in E^+$ ($x_\alpha\xrightarrow{uo(resp. un, uaw)}x$, for short). For further properties about these convergences, see [5-8] for details.

Dunford-Pettis operator has been studied in functional analysis literature. Recently, $uaw$-Dunford?Pettis
operators were introduced and investigated in \cite{OGZ:18,LC:19}. In this paper, we introduce and study the operator which is called \emph{unbounded Dunford-Pettis operator}, and research the relationships with other operators.

A norm on a Banach lattice $E$ is called \emph{order continuous} if $\Vert x_\alpha\Vert\rightarrow0$ for $x_\alpha\downarrow0$. $e\in E_+$ is called an \emph{atom} of the Riesz space $E$ if the principal ideal $E_e$ is one-dimensional. $E$ is called an atomic Banach lattice if it is the band generated by its atoms. A Banach lattice $E$ is said to be $AL$-space if $\Vert x+y\Vert=\Vert x\Vert+\Vert y\Vert$ holds for all $x,y\in E_+$. A Banach lattice $E$ is said to be $AM$-space if $\Vert x+y\Vert=\max\{\Vert x\Vert,\Vert y\Vert\}$ holds for all $x,y\in E_+$. A vector $e>0$ in a Riesz space $E$ is said \emph{strong order unit} whenever the ideal generated by $e$ is $E$.

Recall that an operator $T$ from a Banach space $X$ to a Banach space $Y$ is \emph{Dunford?
	Pettis} if it maps weakly null sequences of $X$ to norm null sequences of $Y$ and is \emph{weak
	Dunford?Pettis} if $f_n(T (x_n)) \rightarrow 0$ for any $w$-null sequence $(x_n)$ in $X$ and any
$w$-null sequence $( f_n)$ in $Y$. An operator $T$ from a Banach lattice $E$ into a Banach space $X$ is said to be an unbounded absolute weak Dunford-Pettis (uaw-Dunford-Pettis, for short)
if for every norm bounded sequence $(x_{n})$ in $E, x_{n}\xrightarrow{uaw}0$ implies $T x_{n}\rightarrow0$ \cite{OGZ:18,LC:19}.

For undefined terminology, notations and basic theory of Riesz spaces and Banach lattices, we refer to [1-4].

\bigskip

\textbf{2. General properties }
\bigskip

\textbf{Definition 2.1} 
An continuous operator $T$ from a Banach lattice $E$ into a Banach lattice $F$ is called \emph{($\sigma$)-unbounded Dunford-Pettis operator} (($\sigma$)-$uDP$, for short) whenever $Tx_\alpha\xrightarrow{un}0(Tx_n\xrightarrow{un}0)$ in $F$ for any bounded $x_{\alpha}\xrightarrow{uaw}0(x_{n}\xrightarrow{uaw}0)$ in $E$. 

An continuous operator $T$ from a Banach lattice $E$ into a Banach lattice $F$ is called \emph{($\sigma$)-unbounded absolute weak continuous operator} (($\sigma$)-$uaw$, for short) whenever $Tx_\alpha\xrightarrow{uaw}0(Tx_n\xrightarrow{uaw}0)$ in $F$ for any bounded $x_{\alpha}\xrightarrow{uaw}0(x_{n}\xrightarrow{uaw}0)$ in $E$.

The collection of all unbounded Dunford-Pettis
operators of $L(E, F)$ will be denoted by $L_{uDP}(E, F)$. That is,
$$L_{uDP}(E, F)=\{T \in L(E, F) : T \text { is uDP operator }\}$$
respectively,
$$L_{\sigma uDP}(E, F)=\{T \in L(E, F) : T \text { is $\sigma$-uDP operator }\}$$
$$L_{uaw}(E, F)=\{T \in L(E, F) : T \text { is $uaw$-continuous operator }\}$$
$$L_{\sigma uaw}(E, F)=\{T \in L(E, F) : T \text { is $\sigma$-$uaw$-continuous operator }\}$$
$$L_{uDP}(E, F)\subset L_{\sigma uDP}(E, F),L_{uaw}(E, F)\subset L_{\sigma uaw}(E, F),$$
$$L_{uDP}(E, F)\subset L_{uaw}(E, F),L_{\sigma uDP}(E, F)\subset L_{\sigma uaw}(E, F).$$

It is clear that a ($\sigma$)-unbounded Dunford-Pettis operator form Banach lattice $E$ into Banach lattice $F$ maps norm bounded relatively (sequentially) $uaw$-compact subset of $E$ to relatively (sequentially) $un$-compact subset of $F$ and a ($\sigma$)-unbounded absolute weak continuous operator form Banach lattice $E$ into Banach lattice $F$ maps norm bounded relatively (sequentially) $uaw$-compact subset of $E$ to relatively (sequentially) $uaw$-compact subset of $F$.

According to $|P_Bx_\alpha-P_Bx|\wedge u\leq|x_\alpha-x|\wedge u$, we have the following result.

\bigskip

\textbf{Lemma 2.2} 
Let $B$ be a projection band and $P$ the corresponding band projection, then $P$ is $uaw$-continuous.

In general, uDP operator does not necessarily exist. For example, $L_p[0,1]$ does not have $uaw$-continuous functional. Next, we show that there exists uDP operator on atomic Banach lattices.

Recall that if $a$ is an atom, then $f_a$ stands for the corresponding \emph{coordinate functional}.

\bigskip

\textbf{Proposition 2.3}
Let $E$ be atomic Banach lattice, $F$ be Banach lattice, if $L(E,F)\ne\{0\}$, then $L_{uDP}(E,F)\ne\{0\}$, moreover $L_{uaw}(E,F)\ne\{0\}$.

\bigskip

\textbf{Proof } 
For any atom $a$ of $E$, let $P_a=f_a(x_\alpha)a$ is a band projection and $P_a$ is $uaw$-continuous, so it is $uaw$-norm-continuous. Considering $T\circ P_e:E\rightarrow F$, it is clear that is uDP operator, moreover it is $uaw$-continuous.

\bigskip
Every unbounded absolute weak Dunford-Pettis operator is $\sigma$-uDP. If the $E^{'}$ has order continuous norm, by \cite[proposition~1]{OGZ:18}, Dunford-Pettis operator is $uaw$-DP, moreover it is $\sigma$-uDP. But the converse is not true in general. Considering the identity operators $I:l_2\rightarrow l_2$, since it is order continuous Banach lattices, by \cite[Theorem 4]{Z:16}, it is uDP operator, but it is not Dunford-Pettis and $uaw$-Dunford-Pettis. 

Dunford-Pettis operator is not necessarily uDP operator. Let $T:l_{1} \rightarrow R$ be defined by $T((x_{n}))=\sum_{n=1}^{\infty}x_{n}$ for every $(x_{n})\in l_{1}$. It is clear that $T$ is compact, hence it is Dunford-Pettis. It follows by considering the standard basis of $l_1$ that $T$ can not be a unbounded Dunford-Pettis operator.	

It's not hard to see that the equivalent characterization of order continuous.
\bigskip

\textbf{Corollary 2.4}
Let $E$ be Banach lattice, then $E$ is order continuous iff $L_{uDP}(E)=L_{\sigma uDP}(E)=L(E)$.
\bigskip

\textbf{Proposition 2.5}
Let $E$, $F$ be Banach lattices, $T:E\rightarrow F$ is an onto lattice homomorphism, if $T$ is Dunford-Pettis operator then $T$ is $\sigma$-uDP operator.	
\bigskip

\textbf{Proof } 
Assume $x_n\xrightarrow{uaw}0$. Since $T$ is onto homomorphism, then for each $u\in E_+$, we have $v\in F_+$ such that $Tu=v$. Thus,
$$T(|x_n|\wedge u)=T|x_n|\wedge v=|Tx_\alpha|\wedge v\rightarrow0$$
for all $u\in E_+$. Hence, $T$ is unbounded Dunford-Pettis.	
\bigskip

\textbf{Corollary 2.6}
Let $E$, $F$ be Banach lattices, $T:E\rightarrow F$ is an onto lattice homomorphism, then $T$ is $uaw$-continuous operator, moreover $T$ is $\sigma$-$uaw$-continuous operator.	
\bigskip

\textbf{3. Ordering structure properties}
\bigskip

 Then,we turn to study the space of uDP operators and $uaw$-continuous operators,let's focus on following questions first.
 
 $(1)$ Whether the uDP oeprators and $uaw$-continuous operators have dominated property?
 
 $(2)$ When the uDP operators (resp. $uaw$-continuous operators) has modulus and the modulus is uDP (resp. $uaw$-continuous)?
 
 $(3)$ In which condition, the space of uDP operators (resp. $uaw$-continuous operators) is order closed?
 \bigskip

\textbf{Proposition 3.1}
For two Banach lattice $E$ and $F$ and two positive operators $T$ and $S$ satisfying $0\leq S\leq T:E\rightarrow F$, if $T$ is uDP, then $S$ is uDP.
\bigskip

\textbf{Proof } 
Since $x_\alpha\xrightarrow{uaw}0$ iff $|x_\alpha|\xrightarrow{un}0$ and $T$ is $s$-continuous operator, so $T|x_\alpha|\rightarrow0$. Using $|Sx_\alpha|\leq S|x_\alpha|\leq T|x_\alpha|$, hence $Sx_\alpha\rightarrow0$, therefore $S$ is $s$-continuous operator.
\bigskip
Similarily, we have the following corollary.

\textbf{Corollary 3.2}
For two Banach lattice $E$ and $F$ and two positive operators $T$ and $S$ satisfying $0\leq S\leq T:E\rightarrow F$, if $T$ is $uaw$-continuous, then $S$ is $uaw$-continuous.
\bigskip

\textbf{Theorem 3.3}
Let $E$ and $F$ be Dedekind complete Banach lattices, for order bounded operator $T:E\rightarrow F$, the following are equivalent:

$(1)$ $T$ is uDP operator.

$(2)$ $T^+$ and $T^-$ are uDP operator.

$(3)$ $|T|$ is uDP operators.
\bigskip

\textbf{Proof } 
$(1)\Rightarrow(2)$ Considering $C_T=\{x\in E:|T|(|x|)=0\}^d$, $C_T$ contains no a disjoint sequence. In fact, for a disjoint sequence $(x_n)\subset C_T$, the $|T|(|x_n|)=\sup T[-|x_n|,|x_n|]$ exists because of $T$ is order bounded and $F$  has order continuous norm. We can find that a disjoint sequence $(y_n)$ satisfying $y_n\in [-|x_n|,|x_n|]$. Since $y_n$ is disjoint, by \cite[Lemma~2]{Z:16}, $\frac{y_n}{\Vert Ty_n\Vert}\xrightarrow{uaw}0$, so $T(\frac{y_n}{\Vert Ty_n\Vert})\rightarrow0$, it is contradiction , so $C_T$ contains finite disjoint atoms at most which denoted by $\{a_1,...,a_m\}$, moreover $C_T=span\{a_1,...,a_m\}$, the $N_T= \{x\in E:|T|(|x|)=0\}$ also is a band, those are projection band because of $E$ is Dedekind complete.

For any bounded $0\leq x_\alpha\xrightarrow{uaw}0$, $Tx_\alpha\xrightarrow{un}0$, 
$T^+x_\alpha=sup\{Ty:0\leq y\leq x_\alpha, y\in C_T\}\leq l \max_{1\leq i\leq m}f_{a_{i}}(x_\alpha)\bigvee^m|Ta_i|$,
since coordinate functional is $uaw$-continuous, so $T^+x_\alpha\xrightarrow{un}0$, hence $T^+$ is uDP, respectively, $T^-$ also is.

$(2)\Rightarrow(3)$ is obvious.

$(3)\Rightarrow(1)$ Using $T=T^+-T^-$, since $0\leq T^+\leq|T|$ and $0\leq T^-\leq|T|$, so $T^+$ and $T^-$ are uDP by Proposition xxx, hence $T$ is uDP.
\bigskip

\textbf{Corollary 3.4}
Let $E$ be Dedekind complete Banach lattice, $F$ be order continuous Banach lattice, for order bounded operator $T:E\rightarrow F$, the following are equivalent:

$(1)$ $T$ is $uaw$-continuous operator.

$(2)$ $T^+$ and $T^-$ are $uaw$-continuous operators.

$(3)$ $|T|$ is $uaw$-continuous operator.
\bigskip

Let $L^b_{uDP}(E,F)$ denotes the collection of order bounded strong continuous operators between two Banach lattices, $L^b_{uaw}(E,F)$ respectively.
\bigskip

\textbf{Corollary 3.5}
Let $E$ and $F$ be Dedekind complete Banach lattices, then $L^b_{uDP}(E,F)$ is a ideal of $L^b(E,F)$.	
\bigskip

\textbf{Corollary 3.6}
Let $E$ be Dedekind complete Banach lattice, $F$ be order continuous Banach lattice, then $L^b_{uaw}(E,F)$ is a ideal of $L^b(E,F)$.
\bigskip

According to the above conclusion, the set of all order bounded $uaw$-continuous functionals in $E^{'}$ is an ideal, but it is not a band in general. Considering $uaw$-continuous functional sequence $f_n:l_1\rightarrow R$, that is,  $f_n(x_n)=\sum_{i=1}^{n}x_i$ for every $(x_n)\in E$, and $f(x_n)=\sum_{i=1}^{\infty}x_i$, it is clear that $f_n\uparrow f$ and $f$ is not $uaw$-continuous.
\bigskip

\textbf{Question 3.7}  
In which condition, $L^b_{uDP}(E,F)$ and $L^b_{uaw}(E,F)$ are bands in $L^b(E,F)$?

\bigskip

\textbf{4. Properties on spaces}
\bigskip

\textbf{Theorem 4.1}
Let $E$ and $F$ be Banach lattices. Then the following assertions are equivalent:

$(1)$ Each positive Dunford?Pettis operator $T:E\rightarrow F$ is unbounded Dunford?Pettis.

$(2)$ Each positive compact operator $T:E\rightarrow F$ is unbounded Dunford?Pettis.

$(3)$ One of the following conditions is valid:

\qquad$(i)$ $E^{'}$ is order continuous.

\qquad$(ii)$ $F=\{0\}$.	
\bigskip

\textbf{Proof } 

$(1)\Rightarrow(2)$ Obvious.

$(2)\Rightarrow(3)$ Assume that $E^{'}$ is not order continuous and $F\neq\{0\}$. We show that there exsits a compact operator which is not
unbounded Dunford?Pettis.

Since $E$ is not order continuous, by Theorem 2.4.14 and Proposition 2.3.11 of \cite{MN:91}, $l_1$ is a closed sublattice of $E$ and there exists a positive projection $P:E\rightarrow l_1$. On the other hand, since $F ={0}$, there exists a vector $0<y\in F_+.$ Define the operator $S:l_1\rightarrow F$ as follows:
$$S\left(\lambda_{n}\right)=\left(\sum_{n=1}^{\infty} \lambda_{n}\right) y$$
for each $(\lambda_n)\in l_1$. Obviously, the operator $S$ is well defined. Let
$$T=S \circ P : E \rightarrow \ell_{1} \rightarrow F$$
then $T$ is a compact operator since $S$ is a finite rank operator $(\text { rank is } 1) .$ But $T$ is not an
unbounded Dunford-Pettis operator. Let $\left(e_{n}\right)$ be the canonical basis of $\ell_{1} .$ Obviously, $e_{n}\xrightarrow{uaw}0$ and $T(e_n)=y\xrightarrow{un}y$. Hence, $T$ is not a unbounded Dunford-Pettis operator.

$(3)(i) \Rightarrow(1)$ Follows from \cite[Proposition~1]{OGZ:18}.

$(3)(i i) \Rightarrow(1)$ Obvious.	
\bigskip

\textbf{Proposition 4.2}
Let $E$ be $AM$-space or atomic order continuous Banach lattice, $F$ be a Banach lattice with strong order unit, if $T:E\rightarrow F$ is $\sigma$-uDP operator then $T$ is Dunford-Pettis.

\bigskip

\textbf{Proof } 
For $x_n\xrightarrow{w}0$, since $E$ is $AM$-space or atomic order continuous Banach lattice, by \cite[Theorem~4.31]{AB:06} and \cite[Proposition~2.5.23]{MN:91}, the lattice operations are weakly sequentially continuous, hence $x_n\xrightarrow{uaw}0$ and $Tx_n\xrightarrow{un}0$. Since $F$ has strong order unit, according to \cite[Theorem~2.3]{KMT:16}, $Tx_n\rightarrow0$.	
\bigskip

\textbf{Corollary 4.3}
Let $E$ be an $AM$-space. Then every operator $T$ from $E$ into arbitrary
Banach lattice is $\sigma$-uDP if and only if $Tx_n\xrightarrow{un}0$ for $x_n\xrightarrow{w}0$.
\bigskip

\textbf{Proof } 
$\Rightarrow$ For any $w$-null sequence $(x_n)$, since $E$ is $AM$-space, then the lattice operations in $E$ is weakly continuous, so $x_n\xrightarrow{uaw}0$, and hence $Tx_n\xrightarrow{un}0$.

$\Leftarrow$ For any bounded $uaw$-null sequence $(x_n)$, since $E^{\prime}$ is $AL$-space which is order continuous, so $x_n\xrightarrow{w}0$ by \cite[Theorem~7]{Z:16}, so $Tx_n\xrightarrow{un}0$, hence $T$ is $\sigma$-uDP.
\bigskip

A continuous operator $T:E\rightarrow F$ between two Banach lattices is said \emph{$\sigma$-unbounded norm continuous} ($\sigma$-$un$-continuous, for short), whenever $x_n\xrightarrow{un}0$ implies $Tx_n\xrightarrow{un}0$ \cite{W:19}.
\bigskip

\textbf{Theorem 4.4}
Let $E$ be atomic Banach lattice with order continuous norm, $F$ be Banach lattice with strong order unit and $T:E\rightarrow F$ be $\sigma$-$un$-continuous operator. Then the following are equivalent:

$(1)$ $E$ is reflexive;

$(2)$ $T$ is Dunford-Pettis operator if and only if $T$ is unbounded Dunford-Pettis operator.
\bigskip

\textbf{Proof } 
It is clear that $F$ is lattice isomorphic to $C(K)$ for compact Hausdorff space $K$.

$(1)\Rightarrow(2)$ Since $E$ is reflexive, it is clear that $E$ and $E^{\prime}$ are $KB$-space. If $T$ is Dunford-Pettis, let $x_n\xrightarrow{uaw}0$, by \cite[Theorem~7]{Z:16},	$x_n\xrightarrow{w}0$ and $Tx_n\rightarrow0$, moreover $Tx_n\xrightarrow{un}0$. If $T$ is unbounded Dunford-Pettis, let $x_n\xrightarrow{w}0$, by \cite[Proposition~2.5.23]{MN:91}, $x_n\xrightarrow{uaw}0$ and $Tx_n\xrightarrow{un}0$, since $F$ has strong order unit, hence $Tx_n\rightarrow0$, therefore $T$ is Dunford-Pettis.

$(2)\Rightarrow(1)$ Assume that $E$ is not reflexive, then $E$ contains no lattice copies of $c_0$ or $l_1$ by \cite[Theorem~4.71]{AB:06}. Suppose that $E$ contains a sublattice isomorphic to $l_1$, let $T:l_{1} \rightarrow R$ be defined by $T((x_{n}))=\sum_{n=1}^{\infty}x_{n}$ for every $(x_{n})\in l_{1}$. $T$ is $\sigma$-$un$-continuous and let $(e_n)$ be standard basis of $l_1$, we have $e_n\xrightarrow{uaw}0$ and $Te_n=1$, hence $T$ is not unbounded Dunford-Pettis, however $T$ is Dunford-Pettis operator. Assume that $E$ contains a sublattice isomorphic to $c_0$ and there exists a lattice complemented isomorphism $S:c_0\rightarrow E$ because of $E$ is order continuous. For $x_n\xrightarrow{uaw}0$ in $c_0$, since $c_0$ is order continuous, then $x_n\xrightarrow{un}0$ by \cite[Theorem~4]{Z:16}, moreover $Sx_n\xrightarrow{un}0$, and since $T$ is $\sigma$-$un$-continuous, hence $TSx_n\xrightarrow{un}0$, so $Tx_n\rightarrow0$ by assumption of $F$ has strong order unit. If $T$ is Dunford-Pettis operators, then $T\circ S$ also is Dunford-Pettis operator, since $c_0$ and $C(K)$ are $AM$-space, so it is contradiction, therefore $T$ is not Dunford-Pettis operator.
\bigskip

Let $E$ be Banach space, $F$ a Banach lattice, and we say that $T:E\rightarrow F$ is \emph{sequentially un-compact} if $TB_E$ is relatively sequentially un-compact \cite{KMT:16}. 
\bigskip

\textbf{Proposition 4.5}
Let $E$, $F$, $G$ be Banach lattices, $T:F\rightarrow G$ be $\sigma$-uDP operator and $S:E\rightarrow F$ be sequentially $un$-compact operator, then $TS$ is sequentially $un$-compact operator from $E$ into $G$.	
\bigskip

\textbf{Proof } 
Suppose $(x_{n})$ is a bounded sequence in $E$. There is a subsequence $\left(x_{n}\right)$ such that $(Sx_{n_{k}})$ is $un$-convergent, so $(Sx_{n_{k}})$ is $uaw$-convergent. Since $T$ is unbounded Dunford-Pettis, hence $(TSx_{n_{k}})$ is $un$-convergent, as desired.	
\bigskip

An operator $T:E\rightarrow F$ is $M$-weakly compact if for every norm
bounded disjoint sequence $(x_n)$ has $Tx_n\rightarrow0$. An operator $T:E\rightarrow F$ is $L$-weakly compact if every disjoint
sequence $(y_n)$ in the solid hull of $T(B_E)$ is norm-null \cite{AB:06}.
\bigskip

\textbf{Proposition 4.6}
Let $E$, $F$ be Banach lattices and $F$ has strong order unit, if $T:E\rightarrow F$ is $\sigma$-uDP operator, then $T$ is $M$-weakly compact.
\bigskip

\textbf{Proof } 
Every bounded disjoint sequence is $uaw$-null in $E$, and since $F$ has strong order unit, hence it is $M$-weakly compact. \cite[Theorem~2.3]{KMT:16}.	
\bigskip

\textbf{Corollary 4.7}
Let $E$, $F$ be Banach lattices and $E$ has order continuous norm, if $T:E\rightarrow F$ is positive $M$-weakly compact, then $T$ is $\sigma$-uDP operator.	
\bigskip
It is similar to \cite[Theorem~4]{OGZ:18}, we have the following conclusion.
\bigskip

\textbf{Corollary 4.8}
Let $E$, $F$ be Banach lattices and $F$ has order continuous norm, then every $L$-weakly compact uaw-continuous operator from $E$ into $F$ is $\sigma$-uDP.	
\bigskip

Recall that an operator $T:E\rightarrow F$ from a Banach lattice to a Banach space is $o$-weakly
compact if it carries order intervals to weakly relatively compact sets.
\bigskip

\textbf{Corollary 4.9}
Every order bounded uDP operator $T:E\rightarrow F$ between two Banach lattices is $o$-weakly compact.
\bigskip

\textbf{Proof } 
Every order bounded disjoint sequence $(x_n)$ is $uaw$-null, and since norm convergent and $un$-convergent are agree for order bounded sequence, we have the conclusion by \cite[Theorem~3.4.4]{MN:91}.
\bigskip

\textbf{Corollary 4.10}
The square of a order bounded uDP operator carries order intervals into relatively sequentially $un$-compact sets.
\bigskip

\textbf{Theorem 4.11}
Suppose $E$ is a Banach lattice with strong order unit and $T$ is a positive uDP operator on
it. If positive operator $S$ dominated by $T^2$, then $S^2$ is compact.	
\bigskip

\textbf{Proof } 
By Proposition 2.8 and Corollary 2.11, $T$ is $o$-weakly compact and $M$-weakly compact. According to Corollary 2.12, $T^2$ maps order intervals into relatively sequentially $un$-compact sets. By \cite[Theorem~2.3]{KMT:16} and \cite[Exercise~5.3.13]{AB:06}, we have the conclusion. 	
\bigskip

By \cite[Proposition~8]{OGZ:18}, we have the following result.
\bigskip

\textbf{Proposition 4.12}
Let $T:E\rightarrow F$ be a positive unbounded Dunford-Pettis operator between Banach lattices
with $F$ Dedekind complete. Then the Kantorovich-like extension $S:E\rightarrow F$ defined via
$$S(y)=\sup \left\{T\left(y \wedge y_{n}\right) :\left(y_{n}\right) \subseteq E_{+}, y_{n}\xrightarrow{uaw} 0\right\}$$
for $y\in E_+$ is again unbounded Dunford-Pettis.
\bigskip

Recall that an operator $T:E\rightarrow F$ between two Banach space is \emph{weak Dunford?Pettis} if $f_n(Tx_n)\rightarrow0$ for any $x_n\xrightarrow{w}0$ in $E$ and any $f_n\xrightarrow{w}0$ in $F^{'}$. A Banach space $X$ is said to have the \emph{Dunford?Pettis property} whenever
$x_{n} \stackrel{w}{\rightarrow} 0$ in $X$ and $x_{n}^{\prime} \stackrel{w}{\longrightarrow} 0$ in $X^{\prime}$ imply $x_{n}^{\prime}\left(x_{n}\right) \rightarrow 0$.
\bigskip

\textbf{Theorem 4.13}
Let $E$ be a Dedekind $\sigma$-complete Banach lattice, we have the following conclusions:	

$(1)$ If $E$ is reflexive, then each weak Dunford?Pettis operator from $E$ into $E$ is unbounded Dunford?Pettis.

$(2)$ If each positive weak Dunford?Pettis operator from $E$ into $E$ is $\sigma$-uDP operator, then $E$ is order continuous.
\bigskip

\textbf{Proof } 
$(1)$ For any weak Dunford-Pettis opeartor $T$ and identity operator $I$, since $E$ is reflexive, so $I$ is weakly compact operator, by \cite[Theorem~5.99]{AB:06}, we have $I\circ T$ is Dunford-Pettis operator.

$(2)$ Assume that $E$ is not order continuous. By \cite[Corollary~2.4.3]{MN:91} and the proof of \cite[Theorem~1]{W:96}, $E$ contains a closed sublattice isomorphic to $l_\infty$ and there is a positive projection $P_1:E\rightarrow l_\infty$, let $S:l_\infty\rightarrow E$ be the
canonical injection and $T=S\circ P_1:E\rightarrow l_\infty\rightarrow E$. Since $l_\infty$ has the Dunford?Pettis property, so $T$ is weak Dunford?Pettis operator. Hence, $T$ is $\sigma$-uDP operator, but considering $(e_n)$, which is the unit vectors of $l_\infty$, it is a contradiction.
\bigskip

\textbf{Theorem 4.14}
Let $E$ and $F$ be Banach lattices. If each weak Dunford?Pettis operator is
$\sigma$-uDP operator, then one of the following assertion is valid:

$(1)$ $E^{\prime}$ is order continuous.

$(2)$ $F$ is order continuous.
\bigskip

\textbf{Proof } 
We show that if $E^{\prime}$ is not order continuous, then $F$ is order continuous. Suppose that $E^{\prime}$ is not order continuous, by Proposition 2.3.11 and Theorem 2.4.14 of \cite{MN:91}, $l_1$ is a closed sublattice of $E$ and there exsits a positive projection $P:E\rightarrow l_1$. Let $(y_n)$ be order bounded disjoint sequence in $F$, define the operator $l_1\rightarrow F$ as: $S(\lambda_n)=\sum_{i=1}^{\infty}\lambda_ny_n$ for each $(\lambda_n)$ of $l_1$, it is well defined, and let $T:S\circ P:E\rightarrow l_1\rightarrow F$. Since $l_1$ has the Dunford-Pettis property, so $T$ is weak Dunford-Pettis operator, hence $T$ is $\sigma$-uDP operator. Considering the standard basis $(e_n)$ of $l_1$, $e_n\xrightarrow{uaw}0$ in $l_1$, so $Te_n=y_n\rightarrow0$ because of $(y_n)$ is order bounded, hence $F$ is order continuous by \cite[Theorem~4.14]{AB:06}.

\end{document}